\documentclass[oneside,a4paper,english, 11pt]{amsart}

\input{artmymacros.sty}
\usepackage{enumitem}
\setlist{nolistsep}

\newcommand{\Frac}{\operatorname{Frac}}

\title{A note on finite length for~$\GL_2$}

\author{Andrea Dotto}
\address{Department of Mathematics,
King's College London,
Strand Campus, Strand, London, WC2R 2LS}
\email{andrea.dotto@kcl.ac.uk}

\author{Daniel Le}
\address{Department of Mathematics,
Purdue University,
150 N. University Street, 
West Lafayette, IN 47907-2067}
\email{ledt@purdue.edu}

\author{Bao V.~Le Hung}
\address{Department of Mathematics,
Northwestern University, 
2033 Sheridan Road\\
Evanston, IL 60208, USA}
\email{lhvietbao@googlemail.com}

\author{Stefano Morra}
\address{LAGA, UMR 7539, CNRS, Universit\'e Paris 13 - Sorbonne Paris Cit\'e, 
Universit\'e de Paris 8,
99 avenue Jean Baptiste Cl\'ement,
93430 Villetaneuse,
France }
\email{morra@math.univ-paris13.fr}
\begin{document}

\begin{abstract}
Let~$p>2$ be a prime and $F/\bQ_p$ be a finite extension.
We prove that every admissible smooth mod~$p$ representation of~$\GL_2(F)$, with a central character, and satisfying a Casimir annihilation condition, 
has finite length.
\end{abstract}

\maketitle


\section{Introduction}
Let $F/\bQ_p$ be a finite extension, and let~$\bF$ be a finite field of characteristic~$p$.
The cohomology of locally symmetric spaces associated to forms of $\GL_2$ over totally real fields is a rich source of
admissible smooth $\bF$-representations~$\pi$ of $\GL_2(F)$.
The study of these representations is a central concern in the mod~$p$ Langlands correspondence for~$\GL_2(F)$, 
but many basic questions remain open.
For example, until very recently, it was not known whether the length of these representations was finite, in any case beyond $F = \bQ_p$.

The papers~\cite{BHHMS2, BHHMSFL} introduced a method for proving that~$\pi$ has finite length, which applies when $F/\bQ_{p}$ is 
unramified and various genericity assumptions on~$\pi$ are met.
This method requires detailed knowledge of the characteristic cycle of~$\pi$, which was established as part of the previous work~\cite{BHHMS}.
Generalizations of~\cite{BHHMS} are not known when $F$ is replaced by a ramified extension of $\Q_p$, nor when $\GL_2$ is replaced by $\GL_n$. 

Recent work of Sorgdrager~\cite{Sorgdrager-thesis}, subsequently refined by Zhao~\cite{Zhao}, gave a new approach to some of the results of~\cite{BHHMS}.
While this approach does not give as much control over the characteristic cycle as~\cite{BHHMS}, it works equally well when $F/\bQ_p$ is ramified, and it gives
a direct computation of the GK-dimension of~$\pi$. 
Combined with \cite[Cor.\ 1.5]{Timmins}, this gives a proof that $\pi$ has finite length when $[F:\Q_p]=2$. 

The aim of this note is to record a very short proof that~$\pi$ has finite length for any extension $F/\bQ_p$, as a corollary of the main result of~\cite{Zhao}.
To state it, we introduce some notation:
\begin{enumerate}
    \item $G := \GL_2(F)$.
    \item $\sm \bF[G]$ is the category of smooth $\bF[G]$-modules.
    \item $T \cong \bZ_p^{[F:\bQ_p]}$ is the commutative Lie subgroup of $G$ defined in~\cite{Zhao}.
    \item $\Lambda_T := \bF[\![T]\!] \cong \bF[\![x_1, \ldots, x_{[F:\bQ_p]}]\!]$.
\end{enumerate}
The paper~\cite{Zhao} also defines a uniform compact open subgroup $H \subset \SL_2(F)$, together with a filtration~$\omega$ on its Iwasawa algebra $\Lambda := \bF[\![H]\!]$.
The graded ring $\gr_\omega \Lambda$ is commutative and Noetherian, and has a distinguished homogeneous ideal~$J$ (the Casimir ideal).

\begin{defn}
    We let~$\cD$ be the full subcategory of $\sm \bF[G]$ consisting of representations~$\pi$ such that: 
    \begin{enumerate}
    \item $\pi^\vee$ is finitely generated over~$\Lambda$, and 
    \item for some (equivalently, by \cite[Lem.\ 2.2]{Zhao}, for any) $\Lambda$-good filtration on~$\pi^\vee$, there exists~$N>0$ such that
    $J^N \gr \pi^\vee = 0$.
    \end{enumerate}
\end{defn}

Note that objects of~$\cD$ are cofinite over~$\Lambda$, hence they are admissible representations of~$G$:
indeed, they are already admissible as representations of $\SL_2(F)$.
The following is our main result.

\begin{thm}\label{thm:main}
    Assume~$p > 2$.
    Let~$\pi \in \cD$. 
    The following are equivalent:
    \begin{enumerate}
    \item $\pi$ has finite $\bF[G]$-length.
    \item $\pi$ is finitely generated over $\bF[G]$.
    \item $\pi$ is $Z$-finite {\em (}in the usual sense that the central subalgebra 
    $\bF[F^\times] \subset \bF[G]$ acts on~$\pi$ via a finite-dimensional quotient{\em )}.
    \end{enumerate}
\end{thm}

Since $Z$-finiteness is implied by the existence of a central character, we immediately get the following corollary.

\begin{cor}\label{cor:central character}
    Assume that~$\pi \in \cD$ has a central character. Then~$\pi$ has finite $\bF[G]$-length.
\end{cor}

    The main application of Cor.\ \ref{cor:central character} arises of course from Sorgdrager's result~\cite[Pf.~of Thm.\ 5.3.1]{Sorgdrager-thesis} 
    that Hecke eigenspaces~$\pi[\fm]$ in the $\bF$-completed cohomology of totally definite quaternion
    algebras over totally real fields, or Shimura curves, are objects of~$\cD$ with a central character, 
    provided that the associated Galois representation~$r_{\fm}$ is absolutely irreducible \cite[Thm.~9.20, Rmk.~9.26]{DPS}. 
    We thus obtain a new proof of the finiteness results in~\cite{BHHMS2, BHHMSFL,Timmins}. 
    However, we emphasize that the methods of~\cite{BHHMS2, BHHMSFL} give more information: for example, they are able to show that~$\pi[\fm]$ is irreducible whenever 
    $F/\bQ_p$ is unramified, and
    the local representation $r_{\fm}|_{G_F}$ is irreducible and sufficiently generic.

\begin{rmk}
    Although~\cite[Thm.\ 1.1]{Zhao} applies equally well when~$\GL_2(F)$ is replaced by the group~$D^\times$ of units in a non-split quaternion algebra over~$F$,
    Cor.\ \ref{cor:central character} is not true for this group, with well-known counterexamples (see e.g.\ \cite[Cor.\ 3.2.4]{BD}). 
    Relatedly, Lem.\ \ref{lem:finite implies semisimple} is not true for this group.
\end{rmk}

\begin{rmk}
After circulating this note, Yongquan Hu pointed out to us that Thm.\ \ref{thm:main} fits the general principle that 
holonomic representations should have finite length, cf.\ \cite[Thm.\ 1.3]{Timmins}.
\end{rmk}

\subsection*{Acknowledgments} 
This note was written without the use of AI tools.
AD was supported by a Royal Society University Research Fellowship during the writing of this note. 
SM was partially supported
by the ANR-25-CE40-4664 and the ANR-MSREI OPE-
2025-0369. 

\section{Proof of Thm.\ \ref{thm:main}}

The proof is a direct consequence of
\cite[Thm. 1.1]{Zhao}, which says that every infinite-dimensional $\pi \in \cD$ has the property that $\pi^\vee$ is finitely generated and faithful over~$\Lambda_T$
(a power series ring), and of the following lemmas.

\begin{lemma}\label{lem:Serre subcategory}
    $\cD$ is a Serre subcategory of $\sm \bF[G]$ {\em (}i.e.\ it is closed under subquotients and extensions{\em )}.
\end{lemma}
\begin{proof}
    Let $0 \to \pi_2 \to \pi \to \pi_1 \to 0$ be an exact sequence in~$\sm \bF[G]$, and let $M := \pi^\vee, N_i := \pi_i^\vee$.
    Since~$\Lambda$ is Noetherian, 
    $N_1, N_2$ are finite over~$\Lambda$ if and only if $M$ is finite over~$\Lambda$.

    Assume that~$\pi \in \cD$.
    Choose a good filtration on~$M$, and equip~$N_1, N_2$ with the induced filtrations, making the sequence 
    $0 \to N_1 \to M \to N_2 \to 0$ strict exact.
    Then the sequence
    $0 \to \gr N_1 \to \gr M \to \gr N_2 \to 0$ is exact (by~\cite[Thm.\ I.4.2.4]{HvO}),
    and
    $\gr M$ is annihilated by~$J^i$ for some~$i$. 
    Hence the same is true of~$\gr N_1, \gr N_2$.
    Since the filtrations on $N_1, N_2$ are good (because $\Lambda$ is complete and $\gr_\omega \Lambda$ is Noetherian, compare \cite[Thm.\ I.5.7]{HvO}), 
    we conclude that $\pi_1, \pi_2 \in \cD$.

    Conversely, if $\pi_1, \pi_2 \in \cD$, then~$M$ is finite over~$\Lambda$, and we can
    choose again a good filtration on~$M$, and equip~$N_1, N_2$ with the induced filtrations. 
    Again, these filtrations are good, and the sequence
    $0 \to \gr N_1 \to \gr M \to \gr N_2 \to 0$ is exact.
    It follows from our assumptions that $J^i$ annihilates~$\gr N_1, \gr N_2$ for~$i$ large enough,
    and so~$J^{2i}$ annihilates~$\gr M$. 
    Since the filtration on~$M$ was chosen to be good, we conclude that~$\pi \in \cD$.
    This concludes the proof.
\end{proof}

\begin{lemma}\label{lem:finite implies semisimple}
Assume~$p > 2$.
If~$\pi \in \sm \bF[G]$ has finite $\bF$-dimension, and $\pi$ has a central character, then~$\pi$ is semisimple. 
\end{lemma}
\begin{proof}
    Since $\GL_2(F)/F^\times\SL_2(F)$ has order coprime to~$p$, it suffices to prove that $\pi$ is semisimple as an $\bF[F^\times\SL_2(F)]$-module.
    Since~$F^\times$ acts on~$\pi$ via a character, it suffices to prove that~$\SL_2(F)$ acts trivially on~$\pi$.
    Let $s:= \left (\begin{smallmatrix}
        0 & {-1} \\ 1 & 0
    \end{smallmatrix} \right )$,
    $t:= \left (\begin{smallmatrix}
        \varpi_F & 0 \\ 0 & 1
    \end{smallmatrix} \right )$, and let $U$ be the upper-triangular unipotent subgroup of $\SL_2(F)$.
    Since $\SL_2(F)$ is generated by~$U$ and $sUs^{-1}$, it suffices to prove that $U$ acts trivially on~$\pi$.
    Since $\pi$ is smooth and finite-dimensional, there is a compact open subgroup~$U^\circ\subset U$ that acts trivially on~$\pi$.
    Since $U = \bigcup_{n \geq 0} t^{-n}U^\circ t^n$, we conclude that~$U$ acts trivially on~$\pi$, as desired. 
\end{proof}

\begin{proof}[Proof of Thm.\ \emph{\ref{thm:main}}]
    (1) implies (2): this is immediate.

    (2) implies (3): since $\pi$ is admissible (being an object of~$\cD$) this is a consequence of~\cite[Lem.\ 2.3.5]{emerton-ordI}.

    (3) implies (1): since~$\pi$ is $Z$-finite, it has a finite filtration whose graded pieces have a central character.
    These graded pieces are objects of~$\cD$, by Lem.\ \ref{lem:Serre subcategory}.
    So it suffices to prove that every $\pi \in \cD$ with a central character has finite length. 
Assume now that $\pi \in \cD$ has a central character. 

    Consider the socle filtration $\soc_i \pi$ of $\pi$, defined so that $\soc_i \pi := 0$ for all $i<0$, 
    $\soc_0 \pi$ is the socle $\soc_{\sm \bF[G]} \pi \subset \pi$, and $\soc_i \pi \subset \pi$ is the preimage of $\soc_{\sm \bF[G]} (\pi/\soc_{i-1}\pi)$ 
    for all $i>0$.  
    The graded pieces of this filtration, which we denote by $\gr_i(\pi) := \soc_i(\pi)/\soc_{i-1}(\pi)$, are objects of~$\cD$, 
    by Lem.\ \ref{lem:Serre subcategory}. 
    Hence they are admissible.
    This implies that $\gr_i(\pi)$ has finite length, for all~$i$.
    (In fact, since every nonzero object of $\sm \bF[G]$ has nonzero $K_1$-fixed vectors, 
    every admissible semisimple object of $\sm \bF[G]$ has finite length.)
    It thus suffices to prove that the socle filtration stops, i.e.\ that $\pi = \soc_N \pi$ for some~$N > 0$. 

    Let~$n := \dim_{\Frac(\Lambda_T)}(\pi^\vee \otimes_{\Lambda_T} \Frac(\Lambda_T))$.
    Since $\pi^\vee$ is finite over~$\Lambda_T$ (by~\cite[Thm.\ 1.1]{Zhao}), this number is finite.
    We will prove that $\pi = \soc_{2n}(\pi)$, i.e.,~that $\gr_{2n+1}(\pi) = 0$.

    If $\gr_i(\pi)$ and $\gr_{i+1}(\pi)$ have finite $\bF$-dimension for some $i\geq 0$ then $\gr_{i+1}(\pi) = 0$:
    indeed, the finite-dimensionality would imply that $\soc_{i+1}(\pi)/ \soc_{i-1} (\pi)$ is finite-dimensional, hence semisimple by Lem.\ \ref{lem:finite implies semisimple}.
    So it suffices to prove that the number of indices~$i$ such that $\gr_i(\pi)$ is infinite-dimensional is at most~$n$.

    We have already observed that $\gr_i(\pi)$ is an object of~$\cD$ for all~$i$. 
    Hence the faithfulness part of~\cite[Thm.\ 1.1]{Zhao} implies that $\gr_i(\pi)$ is infinite-dimensional if and only if
    \[\dim_{\Frac(\Lambda_T)}(\gr_i(\pi)^\vee \otimes_{\Lambda_T} \Frac(\Lambda_T)) > 0.\]
    Since 
    \begin{align*}
n&=\dim_{\Frac(\Lambda_T)}(\pi^\vee \otimes_{\Lambda_T} \Frac(\Lambda_T)) \\
&\geq \dim_{\Frac(\Lambda_T)}(\soc_N\pi^\vee \otimes_{\Lambda_T} \Frac(\Lambda_T)) \\
&= \sum_{i=0}^N \dim_{\Frac(\Lambda_T)}(\gr_i(\pi)^\vee \otimes_{\Lambda_T} \Frac(\Lambda_T))\end{align*}
for any nonnegative integer $N$, the number of strictly positive terms in the final sum is at most $n$ for any $N$. 
\end{proof}

\bibliography{Biblio}

\newcommand{\etalchar}[1]{$^{#1}$}
\def\Dbar{\leavevmode\lower.6ex\hbox to 0pt{\hskip-.23ex \accent"16\hss}D}
  \def\cftil#1{\ifmmode\setbox7\hbox{$\accent"5E#1$}\else
  \setbox7\hbox{\accent"5E#1}\penalty 10000\relax\fi\raise 1\ht7
  \hbox{\lower1.15ex\hbox to 1\wd7{\hss\accent"7E\hss}}\penalty 10000
  \hskip-1\wd7\penalty 10000\box7}
  \def\cfudot#1{\ifmmode\setbox7\hbox{$\accent"5E#1$}\else
  \setbox7\hbox{\accent"5E#1}\penalty 10000\relax\fi\raise 1\ht7
  \hbox{\raise.1ex\hbox to 1\wd7{\hss.\hss}}\penalty 10000 \hskip-1\wd7\penalty
  10000\box7} \def\cprime{$'$}
\providecommand{\bysame}{\leavevmode\hbox to3em{\hrulefill}\thinspace}
\providecommand{\MR}{\relax\ifhmode\unskip\space\fi MR }
\providecommand{\MRhref}[2]{%
  \href{http://www.ams.org/mathscinet-getitem?mr=#1}{#2}
}
\providecommand{\href}[2]{#2}
\begin{thebibliography}{BHH{\etalchar{+}}25b}

\bibitem[BD14]{BD}
Christophe Breuil and Fred Diamond, \emph{Formes modulaires de {H}ilbert modulo
  {$p$} et valeurs d'extensions entre caract\`eres galoisiens}, Ann. Sci. \'Ec.
  Norm. Sup\'er. (4) \textbf{47} (2014), no.~5, 905--974. \MR{3294620}

\bibitem[BHH{\etalchar{+}}23]{BHHMS}
Christophe Breuil, Florian Herzig, Yongquan Hu, Stefano Morra, and Benjamin
  Schraen, \emph{Gelfand-{Kirillov} dimension and mod {{\(p\)}} cohomology for
  {{\(\mathrm{GL}_2 \)}}}, Invent. Math. \textbf{234} (2023), no.~1, 1--128
  (English).

\bibitem[BHH{\etalchar{+}}25a]{BHHMS2}
\bysame, \emph{Conjectures and results on modular representations of {${\rm
  GL}_n(K)$} for a {$p$}-adic field {$K$}}, Mem. Amer. Math. Soc. \textbf{315}
  (2025), no.~1598, v+163. \MR{5003478}

\bibitem[BHH{\etalchar{+}}25b]{BHHMSFL}
\bysame, \emph{Finite length for unramified $\mathrm{GL}_2$},
  \url{https://arxiv.org/abs/2501.03644.pdf}, 2025, preprint.

\bibitem[DPS25]{DPS}
Gabriel Dospinescu, Vytautas Pa\v{s}k\=unas, and Benjamin Schraen,
  \emph{Infinitesimal characters in arithmetic families}, Selecta Math. (N.S.)
  \textbf{31} (2025), no.~4, Paper No. 76, 77. \MR{4940362}

\bibitem[Eme10]{emerton-ordI}
M.~Emerton, \emph{Ordinary parts of admissible representations of {$p$}-adic
  reductive groups {I}. {D}efinition and first properties}, Ast\'erisque
  (2010), no.~331, 355--402. \MR{2667882 (2011k:22013)}

\bibitem[LvO96]{HvO}
Huishi Li and Freddy van Oystaeyen, \emph{Zariskian filtrations},
  $K$-Monographs in Mathematics, vol.~2, Kluwer Academic Publishers, Dordrecht,
  1996. \MR{1420862}

\bibitem[Sor26]{Sorgdrager-thesis}
Reinier Sorgdrager, \emph{{G}elfand--{K}irillov dimensions and the $p$-adic
  {L}anglands program},
  \url{https://pure.uva.nl/ws/files/348670575/Thesis.pdf}, 2026, PhD thesis.

\bibitem[Tim26]{Timmins}
James Timmins, \emph{The canonical dimension of modules for {I}wasawa
  algebras}, Selecta Math. (N.S.) \textbf{32} (2026), no.~1, Paper No. 6, 29.
  \MR{5004711}

\bibitem[Zha26]{Zhao}
Jiusi Zhao, \emph{Ramified tori and {G}elfand--{K}irillov dimension of mod~$p$
  representations}, \url{https://arxiv.org/abs/2611.04350.pdf}, 2026, preprint.
  With an appendix by {Y}ongquan {Hu}.

\end{thebibliography}
\bibliographystyle{amsalpha}

\end{document}